\documentclass[12pt]{amsart}
\usepackage{amssymb,latexsym,eufrak,amsmath,amscd, graphicx}

\theoremstyle{plain}
\newtheorem{theorem}{Theorem}[section]
\newtheorem{proposition}[theorem]{Proposition}
\newtheorem{corollary}[theorem]{Corollary}
\newtheorem{lemma}[theorem]{Lemma}

\newtheorem{theoremalpha}{Theorem}

\theoremstyle{definition}
\newtheorem{definition}[theorem]{Definition}
\newtheorem{remark}[theorem]{Remark}
\newtheorem{example}[theorem]{Example}
\newtheorem{definitionalpha}[theoremalpha]{Definition}

\newcommand{\lra}{\longrightarrow}
\newcommand{\noi}{\noindent}
\newcommand{\PP}{\mathbf{P}}

\newcommand{\RR}{\mathbf{R}}
\newcommand{\NN}{\mathbf{N}}

\newcommand{\ZZ}{\mathbf{Z}}
\newcommand{\CC}{\mathbf{C}}
\newcommand{\QQ}{\mathbf{Q}}

\newcommand{\OO}{\mathcal  {O}}

\newcommand{\fra}{\frak{a}}
\newcommand{\frb}{\frak{b}}
\newcommand{\frc}{\frak{c}}

\newcommand{\frj}{\frak{j}}

\newcommand{\eps}{\varepsilon}

\newcommand{\lbl}{\vskip 6pt}

\newcommand{\ba}{\mathbf{\fra}}

\newcommand{\Bigcone}{\textnormal{Big}}

\newcommand{\bull}{_{\bullet}}

\newcommand{\linser}[1]{\vert \, {#1} \, \vert}
\newcommand{\alinser}[1]{\Vert \, {#1} \, \Vert}

\newcommand{\HH}[3]{H^{{#1}} \big( {#2} , {#3} \big) }

\newcommand{\fall}{ \ \ \text{ for all } \ }

 \newcommand{\MI}[1]{\mathcal  {J} ( {#1} ) }

\newcommand{\pr}{\prime}

\newcommand{\Div}{\text{Div}}

\newcommand{\Biggg}{\textnormal{Big}}

\newcommand{\N}{\mathcal{N}}

\DeclareMathOperator{\Pic}{Pic}

\DeclareMathOperator{\vol}{vol}

\DeclareMathOperator{\ee}{e}

\DeclareMathOperator{\ord}{ord}

\DeclareMathOperator{\Arn}{Arn}

\DeclareMathOperator{\Bs}{Bs}

\DeclareMathOperator{\BB}{{\bf B}}
\DeclareMathOperator{\BBasym}{\BB}
\DeclareMathOperator{\BBstab}{\BB_{+}}
\DeclareMathOperator{\BBrest}{{\bf B}_{--}}
\DeclareMathOperator{\BBig}{Big}

\DeclareMathOperator{\Stab}{Stab}

\begin{document}

\title{Asymptotic invariants of base loci}

\author[L.Ein]{Lawrence Ein}
\address{Department of Mathematics \\ University
of Illinois at Chicago, \hfil\break\indent  851
South Morgan Street (M/C 249)\\ Chicago, IL
60607-7045, USA}
\email{ein@math.uic.edu}

\author[R. Lazarsfeld]{Robert Lazarsfeld} 
\address{Department of Mathematics
\\ University of Michigan \\ Ann Arbor, MI
48109, USA}
\email{rlaz@math.lsa.umich.edu}

\author[M. Musta\c{t}\v{a}]{Mircea Musta\c{t}\v{a}}
\address{Department of Mathematics \\ University of Michigan \\  
Ann Arbor, MI 48109,  USA}
\email{mmustata@umich.edu}

\author[M. Nakamaye]{Michael Nakamaye}
\address{Department of Mathematics and Statistics\\ University of New
  Mexico,
\hfil\break\indent Albuquerque New Mexico 87131,
USA}
\email{nakamaye@math.unm.edu}

\author[M. Popa]{Mihnea Popa}
\address{Department of Mathematics \\
  Harvard University \\ 1 Oxford Street, Cambridge,
  \hfil\break\indent MA 02138,
  USA}
\email{mpopa@math.harvard.edu}

\thanks{2000\,\emph{Mathematics Subject Classification}.
Primary 14C20; Secondary 14B05, 14F17}

\maketitle



\section*{Introduction} 

Let $X$ be a normal complex projective variety, and $D$ a big divisor
on $X$. 
Recall that the \textit{stable base locus}
of $D$ is the Zariski-closed set 
\begin{equation} 
\label{Def.SBL.Intro.Eqn}  \BBasym(D) \ = \ \bigcap_{m > 0} \Bs(mD), 
\end{equation} 
where $\Bs(mD) $ denotes the base locus of the linear
system $\linser{mD}$. This is an interesting and basic invariant, but 
well-known pathologies associated to linear series have discouraged
its study. Recently, however, a couple of results have appeared
suggesting that the picture might be more structured than expected. To
begin with, 
Nakayama \cite{nakayama} attached an asymptotically-defined
multiplicity $\sigma_{\Gamma}(D)$ to any divisorial component $\Gamma$
of $\BBasym(D)$, and proved  
that $\sigma_{\Gamma}(D)$ varies continuously as $D$ varies over the
cone $\Bigcone(X)_{\RR}\subseteq N^1(X)_{\RR}$ of numerical equivalence 
classes of big divisors on $X$. More recently, the fourth author
showed 
in \cite{nakamaye1} that many pathologies disappear if one perturbs
$D$ slightly by subtracting a small ample divisor. Inspired by this
work, 
the purpose of this paper is to define and explore systematically some
asymptotic invariants that one can attach to base loci of linear
series, and to 
study their variation with $D$. 

We start by specifying the invariants in question.
Let $v$ be a discrete valuation of the function field
$K(X)$ of $X$ and let $R$ be the corresponding discrete
valuation ring. Every effective Cartier divisor $D$ on $X$
determines an ideal in $R$, and we denote by $v(D)$
the order via $v$ of this ideal.

Let  $D$ be a big divisor on $X$ with
$\linser{D} \ne \varnothing$. We wish to quantify how nasty are
the singularities of a general divisor $D^\pr \in
\linser{D}$ along
$v$: we define the
\textit{order of vanishing} of $\linser{D}$ along
$v$ by
\begin{align*}
v(\linser{D}) \ &= \ v(D')      
\end{align*}
where $D'$ is
a general divisor in $\linser{D}$.
If $v$ is the valuation associated to a prime divisor $E$ on $X$,
then $v(\linser{D})$ is the coefficient of $E$ in a general element
$D'$ in $\linser{D}$.

Our focus will be on asymptotic analogues
of these invariants. Let $D$ be a big divisor
on $X$, and $v$ a valuation as above. The
\textit{asymptotic order of vanishing} of $D$ along $v$ is
defined as
\begin{equation}
 v(\alinser{D}):= \lim_{p \to \infty}
\frac{v(\linser{pD})}{p}. 
\end{equation} 
It is easy to see that this limit exists.  By taking $p$ to
be sufficiently divisible  this definition extends naturally to
$\QQ$-divisors such that $v(\parallel D\parallel)$ 
is homogeneous of degree one in $D$. 
The most important example is obtained when $X$ is smooth, considering the 
valuation given by the order of vanishing at the generic point of a subvariety $Z$ -- we 
denote the corresponding invariant by $\ord_Z (\parallel D \parallel)$. 
When $Z$ is a prime divisor $E$ on $X$, this is the invariant $\sigma_E(D)$
introduced and studied by Nakayama. In general these invariants may
be irrational (Example \ref{irrational}). 

Our first result
shows that these quantities vary nicely as
functions of
$D$:
\begin{theoremalpha}
\label{Continuity.Theorem.in.Intro}
 Let $X$ be a
normal projective variety, $v$ a fixed
discrete valuation of the function field of $X$,
 and $D$ a big $\QQ$-divisor on $X$. 
\begin{enumerate} \item[(i).] The asymptotic order of vanishing
  $v(\alinser{D})$ depends only on the
numerical equivalence class of $D$, so it induces a function 
on the set $\BBig(X)_{\QQ}$ of
numerical equivalence 
classes of big $\QQ$-divisors.  
\item[(ii).] This function extends uniquely to a continuous function 
on the set 
$\BBig(X)_{\RR}$ 
of
numerical equivalence classes of big $\RR$-divisors. 
\end{enumerate} 
\end{theoremalpha}

When $v$ is the valuation corresponding to a prime divisor 
$E$ on $X$, this
result is due to Nakayama; some of Nakayama's results were
rediscovered and extended to an analytic setting by Boucksom
\cite{boucksom}.  The theorem  
was also suggested  to us by results of the second author
\cite[Chapter 2.2.C]{positivity} 
concerning continuity of the volume of a big divisor. 

It can happen that the stable base locus $\BBasym(D)$ 
does not depend only on the numerical class of $D$ 
(Example~\ref{non-num}). 
Motivated by the work \cite{nakamaye1} of the fourth author, we
consider instead the following approximations of $\BBasym(D)$.

Let $D$ be a big $\QQ$-divisor on $X$. The stable
base locus $\BBasym(D)$ is defined in the natural way, e.g. by taking
$m$ to be sufficiently divisible in (\ref{Def.SBL.Intro.Eqn}). 
The \textit{augmented base locus} of $D$ is the closed set
$$\BBstab(D):=\bigcap_A\BBasym(D-A),$$
where the intersection is over all ample $\QQ$-divisors $A$.
Similarly, the \textit{restricted base locus} of $D$ is given by
$$\BBrest(D):=\bigcup_A\BBasym(D+A),$$
where the union is over all ample $\QQ$-divisors $A$. 
This is a potentially
countable union of irreducible subvarieties of $X$ (it is not known
whether $\BBrest(D)$ is itself Zariski-closed).
It follows easily from the definition that both $\BBrest(D)$ and $\BBstab(D)$
depend only on the numerical class of $D$. 
Moreover, since the definitions involve perturbations
there is a natural way to define the augmented
and the restricted base loci
of an arbitrary real class $\xi\in\BBig(X)_{\RR}$.

The restricted base locus of a big $\QQ$-divisor $D$ is 
the part of $\BBasym(D)$ which is accounted for by numerical properties
of $D$. For example, $\BBrest(D)$ is empty if and only if $D$ is nef.
\footnote{Note that $\BBrest(D)$ appears also in \cite{bdpp}, 
where it is called the \emph{non-nef} locus of $D$.}
At least when $X$ is smooth, we have the following

\begin{theoremalpha}
\label{Theorem.Pos.Stab.Classes} Let $X$
be a smooth projective variety, $v$ a discrete valuation of the function 
field of $X$, and $Z$ the center of $v$ on $X$. If $\xi$ is in 
$\BBig(X)_{\RR}$, then $v(\parallel\xi\parallel)>0$ if and only if
$Z$ is contained in $\BBrest(\xi)$.
\end{theoremalpha}

It is natural to distinguish the big divisor classes for which
the restricted and the augmented base loci coincide. We call such a divisor
class \textit{stable}. Equivalently, $\xi$ is stable if there is a 
neighborhood
of $\xi$ in ${\rm Big}(X)_{\RR}$ such that $\BBasym(D)$ is constant
on the $\QQ$-divisors $D$ with class in that neighborhood. 
For example, if $\xi$ is nef, then it is stable if and only if it is ample.

The set of stable classes is open and dense in
$\BBig(X)_{\RR}$. Given an irreducible closed subset $Z\subseteq X$,
we denote by $\Stab^Z(X)_{\RR}$ the
set of stable classes $\xi$ such that
$Z$ is an irreducible component of $\BBstab(D)$. 
Suppose that $v$ is a discrete valuation with center $Z$.
We see that Theorems~\ref{Continuity.Theorem.in.Intro}  
and~\ref{Theorem.Pos.Stab.Classes} 
show that if $X$ is smooth, then $v(\parallel\cdot\parallel)$ is positive
on $\Stab^Z(X)_{\RR}$ and $v(\parallel\xi\parallel)$
goes to zero only when the argument $\xi$ approaches the boundary of a
connected component of $\Stab^Z(X)_{\RR}$ and when in addition $Z$ 
``disappears'' from the stable base locus as $\xi$ crosses that 
boundary.  We explain the structure of 
the union of these boundaries,  i.e. the set of unstable classes, in 
Section 3. 

Similar asymptotic functions can be defined starting with other invariants of
singularities instead of valuations. 
For example, if $X$ is smooth we may use the reciprocal
of the log canonical threshold or the Hilbert-Samuel multiplicity (cf. Section 2).
The resulting asymptotic invariants enjoy properties analogous to those of 
$v(\parallel\cdot\parallel)$.

In general, these asymptotic invariants
need not be locally polynomial (Example \ref{irrational}). 
However on varieties whose linear series satisfy sufficiently strong 
finiteness hypotheses, the picture is very simple.  We start with a 
definition.
\begin{definitionalpha}[Finitely generated linear series]\label{fg_varieties}
A normal projective variety $X$  has \textit{finitely generated linear
series} if there exist integral Cartier divisors $D_1,\ldots, D_r$ on $X$
with the properties:
\begin{enumerate}\item[(a).] The classes 
of the $D_i$ are a basis for $N^1(X)_{\RR}$;
\item[(b).] The $\ZZ^r$-graded ring 
\[ \text{Cox}(D_1, \ldots, D_r):=\bigoplus_{m = (m_i) \in \ZZ^r}
\HH{0}{X}{\OO_X(m_1D_1 + \ldots + m_rD_r)}\] is a finitely generated
$\CC$-algebra. 
\end{enumerate}
\end{definitionalpha}
\noi The definition was inspired by the notion of a ``Mori dream
space" introduced by Hu and Keel \cite{HK}: these authors require in
addition that the natural map $\Pic(X)_{\QQ} \lra N^1(X)_{\QQ}$ be an
isomorphism, but this is irrelevant for our purposes. It follows 
from a theorem of Cox that any projective toric variety has finitely generated
linear series, and it is conjectured  (and verified in dimension three) \cite{HK} 
that the same is true for any smooth Fano variety. (For more examples
see \S4.)

\begin{theoremalpha}\label{Mori.Dream.Space.Polyhedral.Intro.Thm} 
If $X$  has finitely generated linear series, then the closed cone 
\[ \overline{\textnormal{Eff}}(X)_{\RR}\  = \
\overline{\BBig(X)_{\RR}} \] 
of pseudoeffective divisors\footnote{Recall that by definition a
  divisor is pseudoeffective if its class lies in the closure of the
  cone of effective divisors, or equivalently in the closure of the
  cone of big divisors.} on $X$ is rational polyhedral.
For every
discrete valuation $v$ of the function field of $X$, the function
$v(\parallel\cdot\parallel)$ 
can be extended by continuity to  ${\overline{\rm
    Eff}}(X)_{\RR}$. Moreover, there is a fan $\Delta$
whose support is the above cone, such that for every $v$, the
function $v(\parallel\cdot\parallel)$ is linear on each of the cones in $\Delta$.
\end{theoremalpha} 
\noi A similar statement holds for the usual volume function on varieties
with finitely generated linear series (cf. Proposition \ref{volume}).

The paper is organized as follows. We start in \S 1 with a discussion
of various base loci and stable divisor classes. The asymptotic
invariants are defined in \S 2 and the continuity
is established in \S 3, where we also discuss the structure of the set
of unstable classes and we give a number of examples. Finally, we
prove Theorem \ref{Mori.Dream.Space.Polyhedral.Intro.Thm} in \S 4. 

The present paper is part of a larger project to explore the asymptotic
properties of linear series on X. See [ELNMP1] for an invitation to
this circle of ideas.

\section{Augmented and restricted base loci}

We consider in this section the augmented and restricted base
loci of a linear system, and the notion of stable divisor classes. The
picture is that the stable base locus of a divisor changes only as the
divisor passes through certain ``unstable" classes. 

We start with some notation. Throughout this section $X$ is a normal
complex projective variety. An integral divisor $D$ on $X$ is an
element of the  group $\Div(X)$ of Cartier
divisors, and as usual we can speak about $\QQ$- or $\RR$-divisors. A
$\QQ$- or 
$\RR$-divisor
$D$ is
\textit{effective} if it is a non-negative linear combination of
effective integral divisors with $\QQ$- or $\RR$-coefficients.
  If $D$ is effective, we denote by
${\rm Supp}(D)$ the union of the irreducible
components which appear in the associated
 Weil divisor.  Numerical equivalence
between $\QQ$- or
$\RR$-divisors will be denoted by $\equiv$. We denote by $N^1(X)_{\QQ}$
and 
$N^1(X)_{\RR}$ the finite dimensional $\QQ$- and $\RR$-vector spaces of
numerical equivalence classes.  One has $N^1(X)_{\RR} = N^1(X)_{\QQ}
\otimes_{\QQ} \RR$, and we fix compatible norms
$\parallel\cdot\parallel$ on these two spaces.

Recall next that the \textit{stable base locus} of an integral
divisor $D$ is defined to be
$$\BBasym(D):=\bigcap_{m\geq 1}\Bs(mD)_{\rm red},$$
considered as a reduced subset of $X$. It is elementary that there exists
$p\geq 1$ such that $\BBasym(D)=\Bs(pD)_{\rm red}$, and that  
\begin{equation} \BBasym(D)\ = \ \BBasym(mD) \fall 
 m\geq 1.  \tag{*} \end{equation}
This allows us to define the stable base locus for any $\QQ$-divisor
$D$: take a positive integer  $k$   such that
$kD$ is integral and put
$\BBasym(D):=\BBasym(kD)$. It follows from
(*) that the definition does
not depend on $k$.

\begin{example}[Non-numerical nature of stable base locus]\label{non-num} 
Let $C$ be an elliptic curve, $A$ a divisor of degree $1$ on $C$
and let $\pi\colon X=\PP({\mathcal O}_C\oplus {\mathcal O}_C(A))\to C$.
For $i=1$, $2$, let $L_i={\mathcal O}_X(1)\otimes\pi^*{\mathcal O}(P_i)$, where 
$P_1$ and $P_2$ are divisors of degree zero on $C$, with $P_1$ torsion 
and $P_2$ non-torsion. It is shown in \cite[Example~10.3.3]{positivity}
that $L_1$ and $L_2$ are numerically equivalent big and nef line bundles
such that $\BBasym(L_1)=\varnothing$ and $\BBasym(L_2)$ is a curve.
\qed
\end{example}

\subsection*{\bf The augmented base
locus.} The previous example points to the fact that the stable base
locus of a divisor is not in general very well behaved.  We introduce here
an upper approximation  of this asymptotic locus which has better formal
properties.  The importance of this ``augmented'' locus, and the fact that
it eliminates pathologies, was systematically put in evidence by the
fourth author in \cite{nakamaye1}, \cite{nakamaye2}.

\begin{definition}[Augmented base locus] 
The \emph{augmented base locus} of an  $\RR$-divisor $D$ on $X$ is the
Zariski-closed set:
$$\BBstab(D):=\bigcap_{D=A+E}{\rm Supp}(E),$$
where the intersection is taken over all
decompositions $D=A+E$, where
$A$ and $E$ are $\RR$-divisors such that  $A$ is
ample and $E$ is effective.
\end{definition}

\noindent To relate this definition with the definition given in the introduction,
we note the following:

\begin{remark}[Alternative construction of augmented
base locus]\label{equiv_def} We remark that in the  definition of
$\BBstab(D)$ one may take the intersection over all decompositions such
that, in addition, $E$ is a $\QQ$-divisor. Furthermore,
$$\BBstab(D)=\bigcap_A\BBasym(D-A),$$
 the intersection being taken over all  ample
$\RR$-divisors $A$ such that $D-A$ is a
$\QQ$-divisor.
In fact,  for the first assertion, note
that if $D=A+E$, with $A$ ample and $E$
effective, then we can find effective
$\QQ$-divisors $E_m$ for $m\in\NN$, such that
$E_m\rightarrow E$ when $m$ goes to infinity
and such that ${\rm Supp}(E_m)={\rm Supp}(E)$
for all $m$. Since $D-E$ is ample, so is
$D-E_m$ for $m\gg 0$, hence we are done. The second assertion follows
 immediately from the first one and the definition of the stable
 base locus. 
\end{remark} 

We observe first that -- unlike the stable base locus itself --
the augmented base 
locus depends only on the numerical equivalence class of a divisor.  

\begin{proposition}\label{invar_num} 
If $D_1$ and $D_2$ are numerically equivalent $\RR$-divisors, then
$\BBstab(D_1)=\BBstab(D_2)$. 
\end{proposition} 
\begin{proof} This follows from the observation that if we have a
  decomposition $D_1=A+E$ as in the definition of $\BBstab(D_1)$, then
  we get a corresponding decomposition $D_2=(A+(D_2-D_1))+E$. 
\end{proof}

The next statement shows that $\BBstab(D)$ coincides with  the
stable base locus $\BBasym(D - A)$ for any sufficiently small
ample divisor $A$ such that $D-A$ is a $\QQ$-divisor.  

\begin{proposition}\label{max_baselocus} 
For every $\RR$-divisor $D$, there is $\eps>0$ such that  
\[ \BBstab(D) \ = \ \BBasym(D - A) \] 
for any ample $A$ with $\Vert A \Vert < \eps$ and such that $D-A$ is a
$\QQ$-divisor. More generally, 
if $D'$ is any $\RR$-divisor with $\Vert D'\Vert
<\eps$ and such that
$D-D'$ is a $\QQ$-divisor, then $\BBasym(D-D')\subseteq
\BBstab(D)$.
\end{proposition} 
\begin{proof} There exist ample $\RR$-divisors $A_1,\ldots,A_r$ such
  that each $D-A_i$ is a  $\QQ$-divisor and so that moreover
  $\BBstab(D)=\bigcap_i\BBasym(D-A_i)$. Choose $\eps>0$ so that
  $A_i-D'$ is ample for every $i$ whenever  $\Vert
D'\Vert<\eps$.
  Writing $D-D'=(D-A_i)+(A_i-D')$, we see that if $D-D'$ is a
  $\QQ$-divisor, then $\BBasym(D-D')\subseteq \BBasym(D-A_i)$ for all
  $i$. This proves the second assertion. The first statement follows
  at once, since for any ample divisor $A$ such that $D-A$ is a
  $\QQ$-divisor, we  
have $\BBstab(D)\subseteq
\BBasym(D-A)$ by Remark~\ref{equiv_def}. 
\end{proof}

\begin{corollary}\label{max_baselocus1} 
If $D$ and $\eps>0$ are as in Proposition~\ref{max_baselocus} and if
$D'$ is an $\RR$-divisor such that $\Vert D'\Vert < \eps$, then
$\BBstab(D-D')\subseteq \BBstab(D)$. If $D'$ is ample, then equality
holds.
\end{corollary} 
\begin{proof} For every $D'$ as above, we apply
  Proposition~\ref{max_baselocus} to $D-D'$ to conclude that if $A'$
  is ample, with $\Vert A'\Vert$ small enough, and such that $D-D'-A'$
  is a $\QQ$-divisor, then we have
  $\BBstab(D-D')=\BBasym(D-D'-A')$. Since $\Vert A'\Vert$ is small, we
  may assume that $\Vert D'+A'\Vert <\eps$, hence
  $\BBasym(D-D'-A')\subseteq \BBstab(D)$. Moreover, this is an
  equality if $D'$ (hence also $D'+A'$) is ample. 
\end{proof}  

\begin{example} The augmented base locus is a proper subset
of $X$  if and only if   $D$ is big. Similarly,  it
follows from Proposition~\ref{max_baselocus} that
$\BBstab(D)=\varnothing$ if and only if $D$ is
ample. \qed
\end{example}

\begin{example} \label{Stab.BL.Scaling}
For any $\RR$-divisor $D$, $\BBstab(D) = \BBstab(cD)$ for any real number
$c>0$. \qed
\end{example}

\begin{example} \label{Stab.BL.Sum}
For any $\RR$-divisors $D_1$ and $D_2$, it can be easily shown that 
\[ \BBstab(D_1 + D_2) \ \subseteq \ \BBstab(D_1) \cup \BBstab(D_2). \]
\end{example} 

\begin{example}[Augmented base locus of nef and big
    divisors]\label{nef_case} 
Assume for the moment that $X$ is non-singular, and let $D$ be a nef
and big divisor on $X$. Define the \textit{null locus}
$\text{Null}(D)$ of $D$ to be the union of all irreducible
subvarieties $V \subseteq X$ of positive dimension with the property
that 
$\big( D^{\dim V} \cdot V \big)
= 0$, i.e. with the property that the
restriction of $D$ to $V$ is not big. Then $\BBstab(D) =
\text{Null}(D)$. This is proved for $\QQ$-divisors in
\cite{nakamaye1}, and in general in \cite{ELMNP2}. \qed 
\end{example} 

\begin{example}[Augmented base loci on surfaces] 
Assume here that $X$ is a smooth surface, and let $D$ be a big divisor
on $X$. Then $D$ has a Zariski decomposition $D = P + N$ (see
\cite{Badescu}) into a nef part $P$ and ``negative" part $N$. Then
$\BBstab(D)$ is the null locus $\text{Null}(P)$ of $P$. To see this,
note that if $D=A+E$, where  $A$ is ample and $E$ is effective, then
$E-N$ is effective. 
Therefore $\BBstab(D)=\BBstab(P)\cup {\rm Supp}(N)$.
Since ${\rm Supp}(E)\subseteq\text{Null}(P)$, we
get $\BBstab(D)=\BBstab(P)=\text{Null}(P)$ by the previous
example. \qed 
\end{example}

\subsection*{\bf The restricted base locus.} 
Proposition \ref{max_baselocus} shows that the augmented base locus
of a divisor is 
the stable base locus of a small negative perturbation of the
divisor. When it 
comes time to discuss the 
behavior of the numerical asymptotic invariants of base loci, it will be helpful to
have an analogous notion involving small positive perturbations:

\begin{definition} If $D$ is an $\RR$-divisor
on $X$, then the
\emph{restricted base locus}  of $D$ is
$$\BBrest(D)=\bigcup_A\BBasym(D+A),$$ where the union
is taken over all ample divisors $A$, such that
$D+A$ is a $\QQ$-divisor.
\end{definition}

\begin{remark}[Warning on restricted base loci]\label{warn-rest}
It is not known whether the restricted base locus of a divisor is Zariski
closed in general. \textit{A priori} $\BBrest(D)$ could consist of a
countable union of subvarieties whose Zariski closure is contained in
$\BBstab(D)$. 
\end{remark}

\begin{lemma}\label{equiv_def2} For every
$\RR$-divisor $D$, one has
$\BBrest(D)=\bigcup_A\BBstab(D+A)$, the union
being taken over all ample $\RR$-divisors $A$. 
\end{lemma}

\begin{proof} It is enough to show that if $A$
is ample, then
$\BBstab(D+A)\subseteq
\BBrest(D)$. If $A_0$ is ample with $0<\Vert
A_0\Vert\ll 1$, such that $D+A-A_0$ is a
$\QQ$-divisor, then $A-A_0$ is ample, and 
$$\BBstab(D+A)=\BBasym(D+A-A_0)\subseteq \BBrest(D).$$
\end{proof}

\begin{proposition}\label{invar_num1}
\begin{enumerate}
\item[(i).] For every $\RR$-divisor $D$ and every
real number $c>0$, we have
$\BBrest(D)=\BBrest(cD)$.
\item[(ii).] If $D_1$ and $D_2$ are numerically
equivalent
$\RR$-divisors, then $\BBrest(D_1)=\BBrest(D_2)$.
\end{enumerate}
\end{proposition}

\begin{proof} Both assertions follow from
Lemma~\ref{equiv_def2}, since we already know the
corresponding assertions for the augmented base locus.
\end{proof}

\begin{example} For every $\RR$-divisor
$D$, we have
$\BBrest(D)\subseteq \BBstab(D)$. If $D$ is a
$\QQ$-divisor, then
$\BBrest(D)\subseteq \BBasym(D)\subseteq\BBstab(D)$. \qed
\end{example}

\begin{example}\label{restricted_locus_on_surfaces}
Let $D$ be a big divisor on a smooth projective surface, with Zariski
decomposition $D = P + N$. Then $\BBrest(D) = \text{Supp}(N)$. 
See Example~\ref{restricted_locus_on_surfaces2} below for a proof.
\qed
\end{example}

\begin{example} Given  an $\RR$-divisor $D$, $\BBrest(D)=\varnothing$
 if and only if $D$ is nef. Similarly, $\BBrest(D)=X$ if and only if
 the class of $D$ in $N^1 (X)_{\RR}$ does not lie in the closure of
 the cone of big classes. \qed 
\end{example}

As we have indicated, it isn't known whether $\BBrest(D)$ is
Zariski-closed in general. However it is at worst a countable union of
closed subvarieties: 
\begin{proposition}\label{countable} If
$\{A_m\}_{m\in\NN}$
 are ample divisors with
$\lim_{m\to\infty}\Vert A_m\Vert=0$,
and such that $D+A_m$ are $\QQ$-divisors, then
$\BBrest(D)=\bigcup_m\BBasym(D+A_m)$. In particular,
$\BBrest(D)$ is a countable union of Zariski closed subsets of $X$. 
\end{proposition} 
\begin{proof} The statement follows since for every ample $A$, such
 that $D+A$ is a $\QQ$-divisor,  $A-A_m$ is ample, for $m\gg 0$.
 Since we can write $D+A=(D+A_m)+(A-A_m)$, we get
 $\BBasym(D+A)\subseteq \BBasym(D+A_m)$. 
\end{proof}

\begin{remark}\label{closed_restricted}  
Suppose that $A_m$ is a sequence of ample $\RR$-divisors, where $A_m
\rightarrow 0$. 
As in Proposition~\ref{countable}, we can show that
$\BBrest(D)=\bigcup_m\BBstab(D+A_m)$. We note that if 
$A_m - A_{m+1}$ is ample then $\BBstab(D+A_m) \subseteq
\BBstab(D+A_{m+1})$. In particular, if $\BBrest(D)$ is closed, then it
is equal to $\BBstab(D+A)$ for all sufficiently small ample
$\RR$-divisors $A$. 
\end{remark} 

\begin{proposition}\label{dense_stable} 
For every $\RR$-divisor $D$, there is an $\eps>0$ such that
$\BBrest(D-A)=\BBstab(D-A)=\BBstab(D)$,  for every ample $A$ with
$\Vert A\Vert < \eps$. 
\end{proposition} 
\begin{proof} Apply Corollary~\ref{max_baselocus1} to $D$ to find
  $\eps>0$ such that for every ample $A$, with $\Vert A\Vert < \eps$,
  we have $\BBstab(D-A)=\BBstab(D)$. For every such $A$, we have 
$$\BBstab(D)=\BBstab(D-\tfrac 12
A)\ \subseteq \ \BBrest(D-A)\
  \subseteq\ \BBstab(D),$$ 
as required. 
\end{proof}

\subsection*{\bf Stable divisors.} We now single out those divisors for
which the various base loci we have considered all coincide. 
\begin{definition} An $\RR$-divisor $D$ on $X$
is called
\emph{stable} if $\BBstab(D)=\BBrest(D)$.
\end{definition}

\begin{remark} Note that as both $\BBrest(D)$ and
$\BBstab(D)$ depend only on the numerical class of
$D$, so does the stability condition. 
\end{remark}

The next statement gives various characterizations of stability for an
$\RR$-divisor.

\begin{proposition}\label{equiv_stable} For an
$\RR$-divisor $D$ on $X$, the following are
equivalent:
\begin{enumerate}
\item[(i).] $D$ is stable.
\item[(ii).] There is an ample $\RR$-divisor $A$ such
that
$\BBstab(D)=\BBstab(D+A)$.
\item[(iii).] There is an $\eps>0$ such that for every
$\RR$-divisor
$D'$ with $\Vert D'\Vert<\eps$, we
have
$\BBstab(D)=\BBstab(D+ D')$.
\item[(iv).] There is an $\eps>0$ such that for every
$\RR$-divisor
$D'$ with $\Vert D'\Vert<\eps$, we
have
$\BBrest(D)=\BBrest(D+ D')$. 
\item[(v).]  There exists a positive number  $\eps>0$ such that  all
  the closed sets $\BBasym(D+ D^\pr)$ coincide whenever $D^\pr$ is an
  $\RR$-divisor with $\Vert
D'\Vert<\eps$  and $D+ D^\pr$ rational. 
\end{enumerate} 
\end{proposition} 
\begin{proof} If $D$ is stable, then in particular $\BBrest(D)$ is closed. 
Remark~\ref{closed_restricted} implies that there is a sufficiently
small ample $\RR$-divisor $A$ such that
$\BBstab(D)=\BBrest(D)=\BBstab(D+A)$. This shows that (i)
$\Rightarrow$ (ii).  

We assume (ii). By Corollary~\ref{max_baselocus1}, there is an
$\eps>0$ such that if $\vert D' \vert < \eps$ then $\BBstab(D+
D^\pr)\subseteq \BBstab(D)$. 
We may also assume that $A- D^\pr$ is ample, so $\BBstab(D+A)\subseteq
\BBstab(D+ D^\pr)$. We see that (ii)$\Rightarrow$(iii). 

We assume (iii). Suppose that $\vert D^\pr \vert < \eps$. Note that 
$\BBrest(D+ D^\pr)= \bigcup_A\BBstab(D+ D^\pr +A)$, where $A$ is ample
and $\parallel A \parallel < \eps - \parallel D^\pr \parallel$.
By hypothesis, 
$\BBstab(D+ D^\pr +A)=\BBstab(D)$, and therefore $\BBrest(D+ D^\pr) =
\BBstab(D) =\BBrest(D)$, hence (iv) holds.  

Assume (iv). We choose $\eps$ as in (iv) and such that it satisfies
Proposition~\ref{dense_stable}. 
Assume that $\parallel D'\parallel <\eps$ and that $D+D'$ is
rational. It follows that 
$$\BBrest(D) = \BBrest(D+ D^\pr) \subseteq
\BBasym(D+ D^\pr)
\subseteq \BBstab(D+ D^\pr) \subseteq \BBstab(D).$$
Note that if we take $D'$ such that $-D'$ is ample, then the hypothesis
gives $\BBrest(D)=\BBstab(D)$.
  We deduce that
$\BBasym(D+ D')= \BBstab(D)$ for all $D'$ such that $\vert D'
\vert < \eps$ and $D+ D'$ is rational. 
Hence (iv) impies (v).  

Assume (v). For any sufficiently small ample divisor $A$ such that
$D-A$ is rational,  Proposition \ref{max_baselocus} implies that 
$\BBasym(D-A)=\BBstab(D)$. Now for a sufficiently small ample 
divisor $A'$ such that $D+A'$ is rational, we have that
$$\BBasym(D-A)=\BBasym(D+A')\subseteq \BBrest(D).$$ 
We conclude that $\BBstab(D) = \BBrest(D)$, hence $D$ is stable. 
\end{proof}

\begin{remark} If the class of $D$ is not in the closure of the big
  cone, then $D$ is trivially stable, as $\BBrest(D)=\BBstab(D)=X$. On
  the other hand, if the class of $D$ is in the boundary of this cone
  (so that $D$ is not big), then $D$ is not stable, because
  $\BBstab(D)=X$, but $\BBrest(D)\neq X$. 
\end{remark}

\begin{proposition}\label{density} The set of stable divisor classes
  is open and dense in $N^1 (X)_{\RR}$. In fact, for every $D$ there
  is $\eps >0$ such that if $A$ is ample, and $\Vert A\Vert<\eps$,
  then $D-A$ is stable. 
\end{proposition}

\begin{proof} The set of stable classes is open, as  the condition in
  Proposition~\ref{equiv_stable}(v) is an open condition. To show that
  it is dense, it is enough to prove the last assertion. This follows
  from Proposition~\ref{dense_stable}. 
\end{proof} 

\begin{figure} \label{Stable.Classes.Fig} 
\vskip 5pt
\includegraphics[scale=.8]{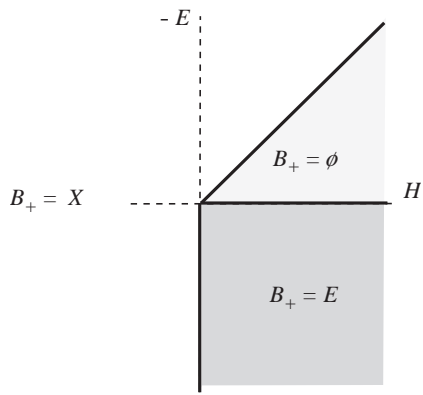} \hskip  .7
in \vskip 5pt \caption{Stable classes on blow-up of $\mathbf{P}^n$} 
\end{figure} 

\begin{example} Let $X = \text{Bl}_{P}(\PP^n) $ be the blowing-up of
  $\PP^n$ at a point $P$. 
Write $H$ and $E$ respectively for the pullback of a hyperplane and
the exceptional divisor. 
For $x , y \in \RR$ consider the $\RR$-divisor $D_{x,y} =  xH - y
E$. Identifying $N^1(X)_{\RR}$ with the $xy$-plane in the evident way,
the set of unstable classes consists of three rays: the negative
$y$-axis, the positive $x$-axis, and the ray of slope $= 1$ in the
first quadrant. The corresponding augmented base loci are indicated
in Figure \ref{Stable.Classes.Fig} . \qed 
\end{example} 

\begin{example} Suppose $X$ is a smooth surface and $D$ a big divisor with
  Zariski decomposition $D = P + N$.  Then $D$ is stable if and only
  if $\text{Null}(P)=\text{Supp}(N)$. \qed 
\end{example} 

\begin{example} For a big $\QQ$-divisor $D$, one introduces in
  \cite{nakamaye2} an asymptotic version of the Seshadri constant at a
  point $x\in X$, denoted by $\eps_m(x,D)$. This invariant describes
  the augmented base locus of $D$, namely $x\in\BBstab(D)$ if and
  only if $\eps_m(x,D)=0$. The main result of \cite{nakamaye2} is a
  continuity statement with respect to $D$ for this invariant. \qed 
\end{example}

\section{Asymptotic numerical invariants} 

In this section we define the asymptotic numerical invariants with
which we shall be concerned. Let $X$ be a normal projective variety.
We fix a discrete valuation $v$ of the function field of $X$,
let $R$ be the corresponding DVR
and $Z$ the center of $v$ on $X$.
Given a big integral divisor $D$, we denote by 
$\fra_p$ the image in $R$
of the 
ideal $\frb(|pD|)$ defining the base locus of $\linser{pD}$. These
ideals form a 
graded system of ideals in the sense of \cite{els}, i.e.
$\fra_p\cdot\fra_{q}\subseteq\fra_{p+q}$ for every $p$ and $q$.
Note that since $D$ is big, we have 
$|pD|\neq\varnothing$, and therefore $\fra_{p}\neq (0)$ for $p\gg 0$.

For every $p$ such that $|pD|\neq\varnothing$, we put 
$v(|pD|)$ for the order $v(\fra_p)$ of the ideal $\fra_p$. Equivalently,
$\nu(|pD|)$ is equal to $v(g)$, where $g$ is an equation 
of a general element in $|pD|$ at 
the generic point of $Z$. 

The convexity property of our invariants 
will be crucial in what follows. A first indication is given 
by the lemma below.
It is an immediate consequence of the fact that the ideals
$\fra_p$ form a graded system of ideals.

\begin{lemma}\label{ineq_for_mult} With the above 
notation, if $p$ and $q$ are such that $|pD|$ and $|qD|$
are nonempty, then
$$v(|(p+q)D|)\leq v(|pD|)+v(|qD|).$$
\end{lemma}

We are now in a position to define our asymptotic invariants.
The existence of the
limit in the following definition follows from Lemma
\ref{ineq_for_mult}. In fact, the limit is equal to the infimum
of the corresponding quantities (see, for example, Lemma~1.4 in \cite{mustata}).
\begin{definition}\label{definition_invariants}
 Given a big integral divisor $D$, set:
$$
v(\parallel D\parallel)= \lim_{p\to \infty}\frac{v(|pD|)}{p}.
$$
This is called the \textit{asymptotic order of vanishing} along $v$.
\end{definition}

\begin{remark}[Rescaling and extension to $\QQ$-divisors]
\label{rescaling} It follows from the definition as a limit that for
any $m \in \NN$:
$$
v(\parallel mD\parallel)=m\cdot v(\parallel D\parallel) \\
$$
In particular,
by clearing denominators we see that 
our invariants are defined 
in the natural way for any big $\QQ$-divisor $D$. 
\end{remark}

\begin{proposition}[Convexity] \label{Convexity.Properties.Invariants}
If $D$ and $E$ are big $\QQ$-divisors on $X$, then
$$v(\parallel D+E\parallel)\leq v(\parallel D\parallel)+v(\parallel E\parallel).$$
\end{proposition}

\begin{proof}
The assertion follows from the fact that if $\fra$, $\frb$ and $\frc$ are
the ideals defining the base loci of the linear systems $|pD|$, $|pE|$
and respectively $|p(D+E)|$, 
then $\fra\cdot\frb\subseteq\frc$.
\end{proof}

\subsection*{\bf Computation via
multiplier ideals.} 
We show now that these asymptotic invariants
can be computed using multiplier ideals.
For the theory of multiplier ideals we refer to
\cite[Part Three]{positivity}.

Note that if $f\colon X'\to X$ is a proper, birational morphism,
with $X'$ normal, then we have an asymptotic order of vanishing
along $v$ defined for big $\QQ$-divisors on $X'$. It is clear that
for a big $\QQ$-divisor $D$ on $X$, we have $v(\parallel D\parallel)
=v(\parallel f^*D\parallel)$. In particular, by taking $f$ such that
$X'$ is smooth, we reduce the computation of the asymptotic order of vanishing
along $v$ to the case of a smooth variety.
In this case, we can make use of multiplier ideals.

Recall that $R$ is the DVR corresponding
to the valuation $v$.
If $D$ is a big integral divisor, we denote by $\frj_p$ the image in $R$ 
of the asymptotic multiplier ideal $\MI{X, \parallel pD \parallel}$.
We show that $v(\parallel D\parallel)$ can be computed using the orders
$v(\frj_p)$ of the ideals $\frj_p$.

The set of ideals $\{\frj_p\}_{p}$ is not a graded sequence anymore.
However, the Subadditivity
Theorem of \cite{del} gives
$\frj_{p+q}\subseteq \frj_p\cdot \frj_q$  for every $p$
and $q$, and hence $v(\frj_{p+q})\geq v(\frj_p)+v(\frj_q)$.
Moreover, if $p<q$  then
$\frj_q\subseteq \frj_p$, so
$v(\frj_q)\geq v(\frj_p)$.
It is easy to deduce from these facts that
\begin{equation}
\lim_{p\to\infty}\frac{v(\frj_p)}{p}=\sup_p\frac{v(\frj_p)}{p}
\end{equation}
(see, for example, Lemma~2.2 in \cite{mustata}).
The above limit is finite: for every $p$ we have
$\fra_p\subseteq\frj_p$ so $v(\frj_p)\leq v(\fra_p)$,
and therefore the above limit is bounded above by $v(\parallel D\parallel)$.
The next proposition shows that in fact we have equality.

\begin{proposition}\label{computed_by_multiplier}
With the above notation, for every big integral divisor $D$
we have
$$v(\parallel D\parallel)=\lim_{p\to\infty}\frac{v(\frj_p)}{p}.$$
\end{proposition}

\begin{proof}
It follows from \cite[Theorem~11.2.21]{positivity}
that there is an effective divisor $E$ on $X$ such that
for every $p\gg 0$, we have
$$\MI{X,\parallel pD\parallel}\otimes{\mathcal O}_X(-E)
\subseteq\frb(|pD|).$$
In particular, there is a nonzero element $u$ in $R$ such that
$u\cdot\frj_p\subseteq\fra_p$ for all $p\gg 0$. 
Therefore $v(\fra_p)\leq v(\frb_p)+v(u)$, so dividing by $p$ and
passing to limit gives
$$v(\parallel D\parallel)\leq\lim_{p\to\infty}\frac{v(\frj_p)}{p}.$$
As we have already seen the opposite inequality, this completes the proof.
\end{proof}

\begin{remark}\label{vanishing_for_b} It
follows from the above proposition that
$v(\parallel D\parallel)=0$ if and only if $\frj_p=R$ for every
$p$. By definition, this is the case if and only if the center $Z$ of $v$
is not contained in ${\mathcal Z}(\MI{X,\alinser{pD}}$
for any $p$.
\end{remark}

\begin{corollary}\label{numerical_invar}
If $D$ and $E$ are numerically equivalent big $\QQ$-divisors 
on the normal projective variety $X$, then $v(\parallel D\parallel)
=v(\parallel E\parallel)$ for every discrete valuation $v$ of the
function field of $X$.
\end{corollary}

\begin{proof}
By taking a resolution of singularities, we may assume that
$X$ is smooth. Moreover, we may assume that
$D$ and $E$ are integral divisors.
Since $D$ and $E$ are big and numerically equivalent,
it follows from \cite[Example 11.3.12]{positivity} that
$$\MI{X,\parallel pD\parallel}=\MI{X,\parallel pE\parallel}$$
for every $p$. The assertion now follows from 
Proposition~\ref{computed_by_multiplier}.
\end{proof}

We can give now a proof
of Theorem~\ref{Theorem.Pos.Stab.Classes} from the Introduction
for $\QQ$-divisors.
In fact, in this case we prove the following more precise 

\begin{proposition}\label{char_restricted_locus}
Let $X$ be a smooth projective variety and $D$ a big $\QQ$-divisor on $X$.
If $v$ is a discrete valuation of the function field of $X$ having
center $Z$ on $X$, then the following are equivalent:
\begin{enumerate} 
\item[(i).]  There is a constant $C>0$ such that $v(|pD|)\leq C$
whenever $|pD|$ is nonempty.
\item[(ii).] $v(\parallel D\parallel)=0$.
\item[(iii).] $Z \not\subseteq \BBrest(D)$.  
\end{enumerate}
\end{proposition}

\begin{proof} 
We may assume that $D$ is an integral divisor.
Note that (i) clearly implies (ii).

Suppose now that $v(\parallel D\parallel)=0$ and let us show that
$Z$ is not contained in $\BBrest(D)$. By Remark~\ref{vanishing_for_b},
we see that $\MI{X,\parallel pD\parallel}={\mathcal O}_X$ at the generic
point of $Z$.
Let $A$ be a very ample divisor on $X$, and $G=K_X+(n+1)A$,
where $n$ is the dimension of $X$. 
It follows from \cite[Corollary 11.2.13]{positivity}
that $\MI{X,\parallel pD\parallel}\otimes {\mathcal O}_X(G+pD)$
is globally generated for every $p$ in $\NN$. This shows that
$Z$ is not contained in
the base locus of $|G+pD|$ for every $p$. Since $G$ is ample,
we deduce from Proposition~\ref{countable} that $Z$ is not contained in
$\BBrest(D)$.

We show now (ii)$\Rightarrow$(i). With the above notation, we see that
$Z$ is not contained in the base locus of $|G+pD|$ for every $p$.
On the other hand, since $D$ is big, by Kodaira's Lemma we can find
a positive integer $p_0$ and an integral effective divisor $B$ such that
$p_0 D$ is linearly equivalent to $G+B$.  
For $p \ge p_0$, 
$pD$ is linearly equivalent with $(p-p_0)D +G+B$, so
$v(|pD|)\leq v(|B|)$. The assertion in (i) follows easily. 

In order to prove (iii)$\Rightarrow$(ii) we proceed similarly.
By Kodaira's Lemma we can find a positive integer $p_0$ and integral
divisors $H$ and $B$, with $H$ ample and $B$ effective such that
$p_0D$ is linearly equivalent to $H+B$. For $p\geq p_0$, we have
$pD$ linearly equivalent to $(p-p_0)D+H+B$. Since
$Z$ is not contained in $\BBrest(D)$, it follows that 
$Z$ is not contained in $\BBstab((p-p_0)D+H)$, hence
$$v(\parallel pD\parallel)\leq v(\parallel (p-p_0)D+H\parallel)+
v(\parallel B \parallel)=v(\parallel B \parallel).$$
Hence $v(\parallel D\parallel)\leq v(\parallel B \parallel)/p$ for every $p$,
and therefore $v(\parallel D\parallel)=0$. 
\end{proof}

\begin{definition}
For every irreducible subvariety $Z$ of a normal variety $X$, there 
is a discrete valuation $v$ whose center is $Z$. For example, take
the normalized blow-up of $X$ along $Z$, and $v$
the valuation corresponding to a component of the exceptional
divisor that dominates $Z$. If $X$ is smooth, we get this way
the valuation given by the order of vanishing at the generic point of $Z$.
For a $\QQ$-divisor $D$, this gives the 
\emph{asymptotic order of vanishing of $D$ along $Z$}: 
$${\rm ord}_Z(\parallel D\parallel):=\lim_{p\to\infty}\frac{{\rm ord}_Z(|pD|)}{p},$$
where ${\rm ord}_Z(|pD|)$ is the order of vanishing along $Z$ of a generic divisor
in $|pD|$.  
\end{definition}

\begin{corollary}
If $X$ is smooth, we have the equality of sets
$$\BBrest(D) = \underset{m\in\NN}{\bigcup}\mathcal{Z}(\MI{X,\alinser{mD}}).$$  
\end{corollary}
\begin{proof}
The assertion follows from Propositions \ref{computed_by_multiplier} 
and \ref{char_restricted_locus}, using the fact that every irreducible 
subvariety $Z$ of $X$ is the center of some discrete valuation. 
\end{proof}

\subsection*{Other invariants}
Similar asymptotic invariants can be defined starting from different
invariants of singularities, instead of valuations. For example, if
$X$ is a smooth variety and $Z$ is an irreducible subvariety, we can consider 
either the Arnold multiplicity or the Samuel
multiplicity at the generic point of $Z$. 
If $R={\mathcal O}_{X,Z}$ and $\fra$ is an ideal in
$R$, then the Arnold multiplicity of $\fra$ at the generic point of $Z$ is
${\rm Arn}_Z(\fra):=1/{\rm lct}(\fra)$, where ${\rm lct}(\fra)$
is the log canonical threshold of the pair $({\rm Spec}(R),V(\fra))$. 
If $\fra$ is of finite colength in $R$, then as usual
$${\rm e}_Z(\fra):=\lim_{\ell\to\infty}\frac{{\rm colength}_R(\fra^{\ell})}
{\ell^d/d!}$$
denotes the Samuel multiplicity of $\fra$, where $d$ is the codimension of $Z$.

Given a big divisor $D$, if $\fra_p$ denotes the image of $b(|pD|)$
in $R$, then ${\rm Arn}_Z(|pD|):={\rm Arn}_Z(\fra_p)$ whenever 
$|pD|$ is nonempty. The corresponding asymptotic invariant is
$${\rm Arn}_Z(\parallel D\parallel):=\lim_{p\to\infty}\frac{{\rm Arn}_Z(|pD|)}
{p}.$$
If $Z$ is not properly contained in any irreducible component of 
$\BBstab(D)$, then for $p\gg 0$ we put ${\rm e}_Z(|pD|)^{1/d}:=
{\rm e}_Z(\fra_p)^{1/d}$ and 
$${\rm e}_Z(\parallel D\parallel)^{1/d}
:=\lim_{p\to\infty}\frac{{\rm e}_Z(|pD|)^{1/d}}
{p}.$$

These invariants satisfy analogous properties 
with $v(\parallel\cdot\parallel)$.
In particular, they can be extended in the obvious way to big $\QQ$-divisors
and depend 
only on the numerical class of the divisor. Moreover, they are positive
at $D$ if and only if $Z$ is contained in $\BBrest(D)$. 
The proofs are analogous using 
the fact that for two ideals $\fra$ and $\frb$ in $R$ we have
$${\rm Arn}_Z(\fra\frb)\leq {\rm Arn}_Z(\fra)+{\rm Arn}_Z(\frb)\,\,
{\rm and}\,\,
{\rm e}_Z(\fra\frb)^{1/d}\leq {\rm e}_Z(\fra)^{1/d}+{\rm e}_Z(\frb)^{1/d}$$
(where in the second inequality we assume that $\fra$ and $\frb$ have
finite colength). The inequality for Arnold multiplicities follows directly
from the definition of log canonical thresholds, while the inequality for
Samuel multiplicities is proved in \cite{teissier}.

\section{Asymptotic invariants as functions on the big cone}

In this section we study the variation of the asymptotic invariants of
base loci. In particular, we prove Theorems
\ref{Continuity.Theorem.in.Intro} and~\ref{Theorem.Pos.Stab.Classes}
from the Introduction. We also present
some examples.

\subsection*{\bf A general uniform continuity lemma.} The continuity
statement (Theorem 
\ref{Continuity.Theorem.in.Intro}) follows formally from an elementary
general statement about convex functions on cones, which we formulate
here. 

Consider an open convex cone $C\subset\RR^n$ and suppose we have a
function $f: C\cap \QQ^n \rightarrow \RR_+$. We assume that this
function satisfies the following properties: 
\begin{enumerate}
\item [(i).] ({\sc homogeneity}). $f(q\cdot x)=q\cdot f(x), ~{\rm
  for~all~} q\in \QQ^*_+$ and $x\in C\cap \QQ^n$. 
\item [(ii).] ({\sc convexity}). $f(x+y)\leq f(x) + f(y),  ~{\rm for~
  all~}  x, y\in C\cap \QQ^n$. 
\item [(iii).] ({\sc ``ample" basis}). There exists a basis $a_1,
  \ldots, a_n$ for $\QQ^n$, 
contained in $C$,   such that $f(a_i) = 0$ for all $i$.
\end{enumerate} 

\begin{proposition}\label{unif_Lipschitz} 
Under the assumptions above, the function $f$ satisfies the following
\textnormal{(}locally uniformly Lipschitz-type\textnormal{)} property
on $C\cap \QQ^n$: for every $x\in  C$, there exists a compact
neighborhood $K$ of $x$ contained in $C$, 
and a constant $M_K>0$, such that for all rational points $x_1, x_2\in
K$ 
\begin{equation} 
|f(x_1)-f(x_2)|\ \leq\ M_K\cdot\parallel
 x_1-x_2\parallel. \tag{$\star$} 
\end{equation} 
In particular, $f$ extends uniquely by continuity to a function on all
of $C$ satisfying $(\star)$.  
\end{proposition} 
\begin{proof} Consider a cube $K\subset C$ with rational
  endpoints. With respect to the chosen basis, we can write 
$$K \ = \  [c_1, d_1  ]\times [c_2, d_2]\times \ldots \times [c_n, d_n].$$ 
We work with the norm given by $\parallel\Sigma u_i a_i\parallel  =
{\rm max}_i \{|u_i|\}$. We have to show that there is $M_K\geq 0$ such
that 
$$|f(x_1)-f(x_2)|\leq M_K\cdot\parallel x_1-x_2\parallel,$$ 
for every $x_1$, $x_2\in K\cap\QQ^n$.

Since $K$ is compact, there exists $\delta \in\QQ_+^*$  such that
$x-\Sigma \delta a_i \in C$ for all $x\in K$. We can also assume (by
subdividing $K$ if  necessary), that all the sides of $K$ have length
$<\delta$. 
Take now any rational points $x_1$, $x_2\in K$, write $x_2 -
x_1=\sum_i \lambda_i a_i$, and set $\lambda= \Vert x_2-
x_1\Vert$. Note that we must have $\lambda < \delta$.  We will estimate the difference $|f(x_1)-f(x_2)|$. 

By repeatedly using properties (i) - (iii) we
get 
\begin{align*} f(x_1)\ &= \ f\big(x_2-(x_2 - x_1)\big)\\ &=\ 
f(x_2-\sum_i\lambda a_i)
\\
  &= \  f\big((1-\lambda /\delta)x_2 + \lambda / \delta
(x_2-\sum_i \delta a_i)\big)\\ 
&\leq \  f\big((1-\lambda /\delta)x_2\big) + f\big(\lambda /
\delta (x_2-\sum_i \delta a_i)\big)\\
&= \ (1-\lambda/ \delta)f\big(x_2\big) + (\lambda/
\delta) f\big(x_2-\sum_i
\delta a_i\big)\\ 
&\leq \ f(x_2) + 
\frac{f(x_2-\sum_i\delta a_i)}{\delta}\cdot
\lambda \\ 
&= \ f(x_2) + 
\frac{f(x_2-\sum_i\delta a_i)}{\delta}\cdot
\Vert x_2 - x_1 \Vert
\end{align*}
On the other hand,  $f(x_2-\sum_i\delta a_i)$
can be bounded  uniformly. Indeed, since
$x_2\in K$, we have that 
$x_2-\sum_i c_i a_i$ is a positive  combination
of the $a_i$'s, so it belongs to $C$ and
$f(x_2-\sum_i c_i a_i)=0$.  Thus we get that 
$$f\big(x_2-\sum_i \delta a_i\big)\ \leq\
f\big(\sum_i(c_i-\delta)a_i\big).$$ If we take
$M_K=f\big(\sum_i(c_i-\delta)a_i\big)/\delta$, it follows that
$$|f(x_1)-f(x_2)|\leq M_K\cdot\parallel
x_1-x_2\parallel$$ 
for every $x_1$, $x_2\in
K\cap\QQ^n$, as required.
\end{proof}

\subsection*{\bf Proof of Theorem \ref{Continuity.Theorem.in.Intro} 
and Theorem \ref{Theorem.Pos.Stab.Classes}.}
We now explain  how Proposition \ref{unif_Lipschitz}
applies to complete  the proofs  of
these two results  from the Introduction. 
Let $v$ be a discrete valuation of the function field of $X$. The 
dependence on the numerical equivalence class in part (i) of
Theorem~\ref{Continuity.Theorem.in.Intro},  
as well as Theorem~\ref{Theorem.Pos.Stab.Classes},  have already been  
proved for $\QQ$-classes in Corollary~\ref{numerical_invar} 
and Proposition~\ref{char_restricted_locus} respectively.

The cone $C$ will be the cone of  big divisors
$\BBig(X)_{\RR}\subset N^1 (X)_{\RR}$.
The fact that the three
properties required in the Proposition are
satisfied for $v(\parallel\cdot\parallel)$  has
already been checked in previous sections:

\begin{enumerate}
\item[(i).] On rational classes, $v(\parallel\cdot\parallel)$
is homogeneous of degree one thanks to  Remark~\ref{rescaling}.

\item[(ii).] The convexity property was noted in Proposition
\ref{Convexity.Properties.Invariants}.  

\item[(iii).]  This follows from the general fact that
one can choose an ample basis for the 
N\'{e}ron-Severi space and the obvious
fact that $v(\parallel A\parallel)=0$ if $A$ is an ample $\QQ$-divisor.
\end{enumerate}

\noi Note that the proof implies the slightly stronger statement that 
the three invariants extend to locally uniformly continuous functions
on the real big cone.

It remains only to show that 
Theorem~\ref{Theorem.Pos.Stab.Classes} holds
for arbitrary big $\RR$-classes $\xi$.  
Suppose first that
$Z\not\subseteq\BBrest(\xi)$, so that for every
ample class
$\alpha$ with $\xi+\alpha$ rational, we have
$Z\not\subseteq\BBrest(\xi+\alpha)$. Corollary~\ref{char_restricted_locus}
gives $v(\parallel \xi + \alpha \parallel)=0$. Letting $\alpha$ go
to $0$, and using continuity, we get $v(\parallel\xi\parallel)=0$. On
the other hand, suppose that $Z\subseteq\BBrest(\xi)$. It follows from
Proposition~\ref{countable} that there is an ample class $\alpha$ such
that $\xi+\alpha$ is rational and
$Z\subseteq\BBrest(\xi+\alpha)$. Therefore
Corollary~\ref{char_restricted_locus} gives
$v(\parallel \xi \parallel)\geq v(\parallel \xi + \alpha \parallel)>0$. 
This completes the proof of Theorems~\ref{Continuity.Theorem.in.Intro}
and~\ref{Theorem.Pos.Stab.Classes}.

\begin{remark}\label{other_invariants2}
If $X$ is a smooth variety and $Z$ is an irreducible subvariety of
codimension $d$,  Proposition~\ref{unif_Lipschitz}
applies for $f={\rm Arn}_Z$, so we get analogues of 
Theorems~\ref{Continuity.Theorem.in.Intro} and~\ref{Theorem.Pos.Stab.Classes}
for ${\rm Arn}_Z(\parallel\cdot\parallel)$. The same applies for 
${\rm e}_Z(\parallel\cdot\parallel)$ with one change: the domain on which it is
defined is ${\rm Big}^Z(X)_{\RR}$, consisting of classes of big $\RR$-divisors
$D$ such that $Z$ is not a proper subset of an irreducible component of
$\BBstab(D)$.
\end{remark}

\lbl 

\noi \textbf{Examples and complements.} We next give some examples and
further information about our invariants. We start with an alternative
computation of the order along a valuation
for a real class. If $v$ is a discrete valuation
of the function field of $X$ and $D$ is an effective divisor on $X$,
then we define $v(D)$ as the order of an equation of $D$ in the local ring 
$R$ of $v$. This extends by linearity to the case of an $\RR$-divisor $D$.

\begin{lemma}\label{alternative_description} 
If $\alpha\in N^1(X)_{\RR}$ is big, then 
\begin{equation}
\label{alternative_formula} v(\parallel\alpha\parallel)=\inf_D v(D), 
\end{equation} 
where the minimum is over all effective $\RR$-divisors $D$ with
numerical class $\alpha$. 
\end{lemma}

\begin{proof} Let us temporarily denote by $v'(\parallel\alpha\parallel)$ the
  infimum in (\ref{alternative_formula}). 
It is easy to check from the definition that $v'$ satisfies
properties (i), (ii) and (iii) in
Proposition~\ref{unif_Lipschitz}. Hence $v'$ is continuous, and
it is enough to show that 
$v'(\parallel\alpha\parallel)=
v(\parallel\alpha\parallel)$ when $\alpha$ is the class of a big
integral divisor $E$. The inequality
$v'(\parallel\alpha\parallel)
\leq v(\parallel\alpha\parallel)$ follows from the definition of the
two functions. 

For the reverse inequality, suppose that $D$ is an effective
$\RR$-divisor, numerically equivalent to $E$. We have to check that
$v(\parallel D\parallel) \leq v(D)$. This is clearly true if
$D$ is a $\QQ$-divisor. In the general case, it is enough to vary the
coefficients of the components of $D$ to get a sequence of effective
$\QQ$-divisors with limit $D$. Taking the limit, we get the desired
inequality.
\end{proof}

Recall that if $X$ is smooth and $Z\subset X$ is an irreducible subvariety,
we denote by ${\rm ord}_Z$ the valuation given by the order of vanishing
at the generic point of $Z$.

\begin{example}\label{restricted_locus_on_surfaces2}
We check the assertion in Example~\ref{restricted_locus_on_surfaces}. 
Let $X$ be a smooth projective surface, and $D$ a big $\RR$-divisor 
with Zariski decomposition $D=P+N$. We prove that $\BBrest(D)={\rm
  Supp}(N)$. If $A$ is ample, then $P+A$ is ample, hence
$\BBstab(D+A)\subseteq {\rm Supp}(N)$. This shows that
$\BBrest(D)\subseteq {\rm Supp}(N)$.

For the reverse inclusion, we use the previous lemma.
If $E$ is an effective $\RR$-divisor numerically equivalent
with $D$, then $E-N$ is effective, so $\ord_Z(\parallel D\parallel)
\geq\ord_Z(N)$ for every $Z$. If $Z$ is a component of $N$, we deduce from
Theorem~\ref{Theorem.Pos.Stab.Classes} that $Z$ is contained in
$\BBrest(D)$.

A similar argument, based on Lemma~\ref{alternative_description},
shows that if $C$ is a curve in $X$, then ${\rm ord}_C(\parallel D\parallel)$
is equal to the coefficient of $C$ in $N$. 
\qed
\end{example}

\begin{example}\label{two_point_blowup}
Let $X={\rm Bl}_{\{P,Q\}}({\mathbf P}^n)$ be the blowing-up of ${\mathbf P}^n$
at two points $P$ and $Q$. 
We assume $n\geq 2$. The N\' eron-Severi group of $X$
is generated
by the classes of the exceptional divisors $E_1$ and $E_2$ and by the
pull-back $H$ of a hyperplane in ${\mathbf P}^n$. A line bundle
$L=\alpha H-\beta_1E_1-\beta_2E_2$ is big if and only if 
$$\alpha>\max\{\beta_1,\beta_2,0\}.$$
We describe now the decomposition of the set of stable classes 
into five chambers and the behavior of our asymptotic invariants on
each of these chambers.

The first region is described by $\beta_1<0$ and $\alpha>\beta_2>0$. 
If $L$ is inside this region, then $L$ is stable and
$\BBasym(L)=E_1$. Moreover, we have $\ord_{E_1}(\parallel
L\parallel)=-\beta_1$. A similar behavior holds inside the second
region, described by $\beta_2<0$ 
and $\alpha>\beta_1>0$. The third chamber is given by $\beta_1$,
$\beta_2<0$ and $\alpha>0$. If $L$ belongs to this chamber, we have
$\BBasym(L)=E_1\cup E_2$, and $\ord_{E_1}(\parallel
L\parallel)=-\beta_1$ and 
$\ord_{\parallel E_2 \parallel}(L)=-\beta_2$.

{}From now on we assume that $\beta_1$, $\beta_2>0$. The fourth chamber 
is given by adding the condition $\alpha>\beta_1+\beta_2$. This
chamber gives precisely the ample cone. The last region is given by
the opposite inequality $\alpha<\beta_1+\beta_2$. Every $L$ in this
chamber is stable, and $\BBasym(L)=\ell$, the proper transform of the
line $\overline{PQ}$. In order to compute the invariants associated to
$L$ along $\ell$, we may assume that $P=(1:0:\ldots:0)$ and
$Q=(0:1:0:\ldots:0)$. 
We see that $H^0(X,L)$ is spanned by 
$$\{\prod_{i=0}^nX_i^{p_i}\vert
p_0\leq\alpha-\beta_1,p_1\leq\alpha-\beta_2, \sum_ip_i=\alpha\}.$$ 
Therefore we get coordinates $x_2,\ldots,x_{n}$ at the generic point 
of $\ell$ such that the base locus of $L$ is defined at this point by 
$(\prod_{i=2}^nx_i^{q_i}\vert\sum_iq_i\geq \beta_1+\beta_2-\alpha)$
and  $\ord_{\ell}(\parallel L\parallel)=\beta_1+\beta_2-\alpha$.
Note that we also have 
$\Arn_{\ell}(\parallel L\parallel)=(\beta_1+\beta_2-\alpha)/(n-1)$ and
$\ee_{\ell}(\parallel L\parallel)^{1/(n-1)}=(\beta_1+\beta_2-\alpha)$.

We will see in the next section that for every toric variety 
(or more generally, for every variety with finitely generated linear
series) there is a fan refining the big cone as above, such that on
each of the subcones our asymptotic invariants are polynomial.\qed
\end{example}

\begin{example}\label{irrational}
We give now an example when the asymptotic invariants can take 
irrational values for $\QQ$-divisors. Moreover, we will see that 
in this case the invariants are not locally polynomial. 
The idea of this example is due to Cutkosky \cite{cutkosky}. We follow
the approach in K\"{u}ronya \cite{kuronya} where this is used to give
an example when the volume function is not locally polynomial.

We start by recalling the notation and the definitions from
\cite{kuronya}. Let $S=E\times E$, where $E$ is a general elliptic
curve. If $F_1$ and $F_2$ are fibers of the respective projections,
and if $\Delta$ is the diagonal, then the classes of $F_1$, $F_2$ and
$\Delta$ span $N^1(X)_{\RR}$. If $h$ is an ample class on $S$ and if 
$\alpha\in N^1(X)_{\RR}$, then $\alpha$ is ample (equivalently, it is
big) if and only if $(\alpha^2)>0$ and $(\alpha\cdot h)>0$. We
consider the following ample divisors on $S$: $D=F_1+F_2$ and
$H=3(F_2+\Delta)$.

Let $\pi : X=\PP(\OO_S(D)\oplus\OO_S(-H))\longrightarrow S$ 
be the canonical projection. If $0\leq t\ll 1$, with $t\in\QQ$, we
take  $L_t=\OO(1)+t\cdot\pi^*F_1$, which is big. We consider the
section of $\pi$ induced by the projection
$\OO_S(D)\oplus\OO_S(-H)\longrightarrow\OO_S(-H)$, and denote by $E$
its image. 
We will compute $\ord_E(\parallel D_t\parallel)$. If $k$ is a positive
integer such that $kt\in\NN$, then
$$H^0(X,\OO_X(kD_t))\simeq\oplus_{i+j=k}H^0(S,\OO_S(iD-jH+ktF_1)).$$

An easy computation shows that if 
$$\sigma(t)=\frac{9+5t-\sqrt{49t^2+78t+45}}{18-12t},$$
then $H^0(S,\OO_S(iD-jH+(i+j)tF_1))$ is zero if $j/i>\sigma(t)$ 
and it is non-zero if $j/i<\sigma(t)$. Note also that 
$\OO_X(kD_t-pE)\vert_E\simeq\OO_S(pD-(k-p)H+ktF_1)$. We deduce 
$$|\ord_E(|kD_t|)-\lceil k/(1+\sigma(t))\rceil|\leq 1.$$
This implies $\ord_E(\parallel D_t\parallel)=1/(1+\sigma(t))$. By
taking $t=0$, we get $\ord_E(D_0)\not\in\QQ$. Moreover, it is clear
that $\ord_E$ is not a locally polynomial function in any
neighbourhood of $D_0$. \qed
\end{example}

\begin{example}[Surfaces]\label{case_of_surfaces}
The case of surfaces has been studied recently both from the 
point of view of the volume function and of asymptotic base loci
in \cite{BKS}. We interpret now their results in our framework.

Let $X$ be a smooth projective surface, and let $B\subset X$ be a
collection of curves  having negative definite intersection
form. Consider the (possibly empty) set ${\mathcal S}_B$ consisting of
stable classes $\alpha\in N^1(X)_{\RR}$ with $\BBstab(\alpha)=B$. It
is clear that if it is non-empty then ${\mathcal S}_B$ is an open
cone. Moreover, it is also convex, since if $\alpha_1=P_1+N_1$ and
$\alpha_2=P_2+N_2$ 
are two Zariski decompositions with ${\rm Supp}(N_1)={\rm Supp}(N_2)$, 
then $\alpha_1+\alpha_2=(P_1+P_2)+(N_1+N_2)$ is the Zariski
decomposition of $\alpha_1+\alpha_2$. 

If $E_1,\ldots,E_r$ are the irreducible components of $B$, and if 
$\alpha\in {\mathcal S}_B$ has Zariski decomposition
$\alpha=P+\sum_{j=1}^ra_jE_j$, then $(\alpha\cdot
E_i)=\sum_{j=1}^r(E_i\cdot E_j)a_j$. Therefore the coefficients 
$a_j$ depend linearly on $\alpha$, hence for every curve $C$ on $X$,
the function $\ord_C(\parallel\cdot\parallel)$ 
is linear on ${\mathcal S}_B$, with rational
coefficients. 
 
The closed cones $\overline{{\mathcal S}_B}$ give a cover of
$\overline{\BBig(X)}_{\RR}$
that is locally finite 
inside the big cone. 
Indeed, suppose that $\alpha\in\BBig(X)_{\RR}$. 
It follows from Corollary~\ref{max_baselocus1} that if $\beta$ is in a
suitable open neighbourhood $U$ of $\alpha$ and if $\beta\in {\mathcal S}_B$,
then $B\subseteq \BBstab(\alpha)$. In particular, there are only
finitely many possibilities for $B$. 

We show now that each cone $\overline{{\mathcal S}_B}$ is rational
polyhedral inside the big cone. 
We keep the above notation, and
without any loss of generality, we may assume that the
open subset $U$ is a convex cone. We have seen that $U$ is covered by
finitely many $\overline{{\mathcal S}_{B_i}}$, and
for every curve $C$ in $X$,
we have linear functions $L_i$ such that ${\rm ord}_C(\parallel\cdot\parallel)=L_i$
on $\overline{{\mathcal S}_{B_i}}$. Since the asymptotic order function 
along $C$ is convex,
it follows from general considerations that
${\rm ord}_C(\parallel\cdot\parallel)=\max_iL_i$ on $U$, i.e. 
${\rm ord}_C(\parallel\cdot\parallel)$ is piecewise linear on $U$. 
On the other hand,
${\mathcal S}_B\cap U$ is the set of 
those $\beta\in U$ such that
${\rm ord}_C(\parallel\xi\parallel)=0$ for
$\xi$ in a neighborhood of $\beta$, for every
$C$ in $\BBstab(\alpha)$ but not in $B$, and ${\rm ord}_C(\parallel\beta\parallel)
\neq 0$ for $C$ in $B$. Therefore $\overline{{\mathcal S}_{B}}\cap U$
is the intersection of
$U$ with finitely many half-spaces, which proves our assertion.

\qed
\end{example}

\subsection*{\bf The structure of the unstable locus.}
We discuss the structure of the locus of unstable  classes
inside  the big cone. The picture is similar to that given by a theorem
of  Campana and Peternell (cf.  
\cite[Chapter 1.5]{positivity})   for the
structure of the boundary of the nef cone.
We assume that $X$ is smooth.

We first fix a closed subset 
$Z\subseteq X$.
Let $v_Z$ be a fixed discrete valuation of the function field of $X$
such that $Z$ is the center of $v_Z$ on $X$.
We use the asymptotic order
function $v_Z(\parallel\cdot\parallel)$  to obtain information on the
locus of $Z$-unstable points.  The zero locus
\[  \N_{Z} :=\{\xi\in {\rm Big}(X)_{\RR}\mid
v_Z(\parallel\xi\parallel)=0\}\] is a convex cone
which is closed in
$\BBig(X)_{\RR}$. 
By Theorem~\ref{Theorem.Pos.Stab.Classes} this is the set of big classes
$\xi$ such that $Z$ is not contained in $\BBrest(\xi)$.
We call it the null cone determined by $Z$.
A class $\xi \in \BBig(X)_{\RR}$  is called
$Z$-\emph{unstable} if
$Z\subseteq \BBstab(\xi)$,  but $Z \not
\subseteq \BBrest(\xi)$.
It is easy to see that
the $Z$-unstable classes are precisely the big classes
that lie on the boundary of ${\mathcal N}_{Z}$.

By definition, a class 
$ \xi \in \BBig(X)_{\RR}$ is unstable if and  only if it
is $Z$-unstable for some irreducible  component
$Z\subseteq \BBstab(\xi)$.  Thus $\xi$ is unstable
if and only if it is $Z$-unstable for \emph{some}
subvariety $Z$.   Thus the picture is that we 
have convex null-cones $\N_Z$ in $\BBig(X)_{\RR}$
indexed by  all subvarieties $Z\subseteq X$, and 
\begin{equation}  \textrm{Unstab} (X) \ = \ \underset{Z}{\bigcup}
~\partial\N_Z.\tag{*} \end{equation}  
It follows for example that the set of unstable classes
does not contain isolated rays.
(This is just a general statement
about boundaries of convex cones. Visually, this says that in any section 
of the big cone the unstable locus does not have isolated points.)

Note that the union in $(*)$ can be taken over countably many $Z$.
Indeed, it is enough to consider those $Z$ which are irreducible components 
of augmented base loci (and we may restrict to $\QQ$-divisors by
Proposition~\ref{max_baselocus}). Since $\BBstab(D)$ depends only on
the numerical equivalence class of the $\QQ$-divisor 
$D$, we have to consider only countably many subvarieties.

Since $\partial \N_Z$ is the boundary of a convex
cone in the N\'{e}ron-Severi space,  it has measure
zero and therefore so does $\textrm{Unstab}(X)$.

\begin{remark}
In fact, one can show that 
there is an open 
dense subset $V\subset {\rm Unstab} (X)$ which looks locally like the 
boundary of a unique ${\mathcal N}_Z$: for every
$\xi\in V$, there is an open neighborhood $U(\xi)$ of $\xi$,
and an irreducible closed subset $Z\subseteq X$, such that
${\rm Unstab} (X) \cap U(\xi)=\partial {\mathcal N}_Z\cap U(\xi)$.
\end{remark}

\section{Asymptotic invariants on varieties with finitely generated
  linear series}

Our goal in this section is to prove Theorem
\ref{Mori.Dream.Space.Polyhedral.Intro.Thm} from the Introduction. In
fact, we 
will prove a somewhat stronger local statement.

Let $X$ be a normal projective variety, and fix $r$ integral divisors
$D_1, \ldots, D_r $ on $X$ 
such that some linear combination of the $D_i$ (with rational
coefficients) is big.  Setting $N = \ZZ^r$ and $N_{\RR} = N \otimes
\RR = \RR^r$, 
the choice of the $D_i$ gives linear maps
\[  \phi : N \lra N^1(X) \ \quad , \ \quad \phi_{\RR} : N_{\RR} \lra
N^1(X)_{\RR}. \] 
We denote by $B \subseteq N_{\RR}$ the pull-back
$\phi_{\RR}^{-1}\big(\Biggg(X)_{\RR} \big)$, 
so that $\overline{B}$ is the pull-back of the closure of
$\Biggg(X)_{\RR}$. The main result of this section is: 

\begin{theorem}\label{MDS_local}
Assume that the graded $\CC$-algebra \begin{equation}
  \label{Local.Cox.Ring.Eqn} \textnormal{Cox}(D_1, \ldots, D_r):= 
\bigoplus_{m=(m_i)\in\ZZ^r}\HH{0}{X}{\OO_X(m_1D_1 + \ldots + m_r D_r)}
\end{equation}
is finitely generated. Then $\overline{B}$ is a rational polyhedral
cone and for every discrete valuation $v$ of the function field of $X$,
the pull-back to $B$ of the function $v(\parallel\cdot\parallel)$ can be extended
by continuity to $\overline{B}$. Moreover, there is a fan $\Delta$ with
support $\overline{B}$ such that every $v(\parallel\cdot\parallel)$
is linear on the cones in $\Delta$.
\end{theorem}

Before proving Theorem \ref{MDS_local} we give a few examples of
finitely generated Cox rings.  

\begin{example}  Suppose $N^1(X)_\RR$ has dimension 1.  Let $D$ be any
  ample divisor on $X$.  Then the $\ZZ$--graded ring
$$
{\text{Cox}}(D) = \bigoplus_{m \in {\bf Z}} H^0(X,mD)
$$
is finitely generated since it is isomorphic to the projective
coordinate ring of $X$. Hence $X$ has finitely generated linear
series. \qed
\end{example}

\begin{example}  If $X$ is the projective plane blown up at an
arbitrary number of \emph{collinear} points, 
then it is shown in \cite{EKW} that $X$ has finitely generated linear
series.
\end{example}

\begin{example}
If $X = {\rm Bl}_{p_1, \ldots,p_r} (\PP^n)$, where $n\geq 2$, $r\geq
n+3$, and $p_1, \ldots, p_r$ are distinct points lying on a rational 
normal curve in $\PP^n$, then it is shown in \cite{CT} that $X$ has
finitely generated linear series. (Cf. \emph{loc. cit.} for a few
other examples.)  
\end{example}

To prove Theorem~\ref{MDS_local}, it is convenient to pass to a local statement
involving families of ideals.

\begin{definition}\label{definition_multi-graded} 
Let $V$ be any  variety, and let $S\subseteq \ZZ^r$ be a subsemigroup 
(in most cases we will take $S=\NN^r$ or $S=\ZZ^r$). An
\textit{$S$-graded system of ideals} on $V$ is a collection $\fra\bull
= \{ \fra_m \}_{m \in S}$ 
of ideal sheaves on $V$, with $\fra_0 = \OO_V$, which satisfies 
\[  \fra_m \cdot \fra_{m^\pr}  \subseteq \fra_{m + m^\pr}\]
for all $m , m^\pr \in S$. The \textit{Rees algebra} of $\fra\bull $
is the $S$-graded $\OO_V$-algebra 
\[ \mathit{R}(\fra\bull) \ = \ \bigoplus_{m \in S} \fra_m, \] 
and $\fra\bull$ is \textit{finitely generated} if
$\mathit{R}(\fra\bull)$ is a finitely generated $\OO_V$-algebra. \qed
\end{definition}

\noi For example, starting with divisors $D_1, \ldots, D_r $ on a 
projective variety $X$, we have 
an $\NN^r$-graded sequence of ideals $\frb\bull$
such that $\frb_m$ is the ideal defining the base locus of
$|m_1D_1+\ldots+m_rD_r|$.
If the Cox ring $\text{Cox}(D_1, \ldots, D_r)$
in (\ref{Local.Cox.Ring.Eqn}) is finitely generated, then the
corresponding system $\frb\bull$  of base ideals is likewise finitely
generated.

\begin{remark}[Invariants for graded systems]
\label{regularity_properties}  
It would be very interesting to know what sort of regularity
properties the functions defined by the invariants introduced
in \S2 satisfy. For example, are they piecewise analytic on a 
dense open set in their domains?  
As the reader has probably noticed, these invariants 
can also be defined for an arbitrary graded sequence of ideals.
As a consequence, most of what we have done in the previous sections can be transposed into the abstract setting of $S$-graded systems (in this case the N\'{e}ron-Severi space is replaced by the group generated by $S$). 
One can define analogues of the effective and of the nef cones
in this setting and under mild hypotheses (for example that the system in question
contain a non-empty ``ample'' cone), one can prove the
continuity of the asymptotic invariants in this setting.
See \cite{Wolfe} and \cite{ELMNP1} for more on this point of view.
Work of Wolfe \cite{Wolfe} suggests that in this abstract 
setting one can't generally expect any good behavior other than that 
implied by convexity. One might hope however that this sort of
pathology does not occur in the global geometric setting. 
\end{remark}

For the proof of Theorem \ref{MDS_local}, the main point is to
show that a finitely generated graded system essentially 
is given by products of powers of finitely many ideals. This is the content
of the following Proposition. We fix a lattice $N\simeq\ZZ^r\subset
N_{\RR}=N\otimes_{\ZZ}\RR$ 
and a finitely generated, saturated subsemigroup $S\subseteq N$. 
This means that if $C$ is the convex cone generated by $S$, then 
$C$ is a rational, polyhedral cone, and $S=C\cap N$. We denote
by $\overline{\fra}$ the integral closure of an ideal $\fra$.

\begin{proposition}\label{MDS_algebraic} 
With the above notation, let $\ba_{\bullet}$ be a finitely generated
$S$-graded system of ideals on the variety $V$ (more generally, $V$
can be an arbitrary Noetherian scheme). Then there is a smooth fan
$\Delta$ with support $C$, such that for every smooth refinement $\Delta'$
of $\Delta$ there is a positive integer $d$ with the following
property. For every cone $\sigma\in\Delta'$,  if we denote by
$e_1,\ldots,e_s$ the generators of
$S_\sigma:=\sigma\cap N$,  then 
\begin{equation}\label{inclusion_in_product}
\overline{\ba_{d{\sum_i}p_ie_i}}=\overline{\prod_i
\ba_{de_i}^{p_i}},
\end{equation} 
for every $p=(p_i)\in\NN^s$.  
\end{proposition}

It is clear that it is enough to prove Proposition~\ref{MDS_algebraic} when $X={\rm
  Spec}(R_0)$ is affine. Before giving the proof we need a few
  lemmas. The following one is well known, but we include a proof for
  the benefit of the reader.

\begin{lemma}\label{subcone} 
With $S$ as above, suppose $R=\oplus_{m\in S}R_m$ is an $S$-graded
ring that  is finitely generated as an $R_0$-algebra. If $S'\subseteq
S$ is a (finitely generated, saturated) subsemigroup, and if
$R'=\oplus_{m\in S'}R_m$, then $R'$ is a finitely generated
$R_0$-algebra.
\end{lemma}

\begin{proof} Choose homogeneous generators $x_1,\ldots,x_q$ of $R$ as
  an $R_0$-algebra, and let $m_i=\deg(x_i)$. We get a surjective
  morphism of $R_0$-algebras
$$\Phi : R_0[X_1,\ldots,X_q]\longrightarrow R$$ 
given by $\Phi(X_i)=x_i$. This is homogeneous with respect to the
semigroup homomorphism $\phi : \NN^q\longrightarrow S$ that takes the
$i^{\rm th}$ coordinate vector to $m_i$. 

If $T:=\phi^{-1}(S')$, then $T$ is cut out in $\NN^q$ by finitely many
linear inequalities, hence $T$ is finitely generated by Gordan's
Lemma. If $w=(w_1, \ldots,w_q)\in\NN^q$, we put $X^w$ for the monomial
$\prod_iX_i^{w_i}$. If we choose generators $v^{(1)},\ldots,v^{(p)}$
for $T$, and let $y_i=
\Phi(X^{v^{(i)}})$, then $R'$ is generated
over $R_0$ by $y_1,\ldots,y_p$. For this, it is enough to note that by
the surjectivity of $\Phi$, every homogeneous element in $R'$ is a
linear combination (with coefficients in $R_0$) of images of monomials
with degrees in $T$.
\end{proof}

We will prove Proposition~\ref{MDS_algebraic} by induction on
$\dim(S)$.  The following lemma which covers the case $S=\NN$ is
standard (see \cite{Br}, Chapter III, Section 1, Proposition~2). Note
that in this case we get a stronger statement than in
Proposition~\ref{MDS_algebraic}.

\begin{lemma}\label{dim_one}
If $R_0$ is a Noetherian ring, and if $R=\oplus_{m\in\NN}R_m$ is an
$\NN$-graded, finitely generated $R_0$-algebra, then there is a
positive integer $d$, such that $R_{dm}= R_d^m$, for every
$m\in\NN\setminus\{0\}$.
\end{lemma}

We need one more easy result about cones and semigroups.

\begin{lemma}\label{multiple_in_cones} Let $N\subset N_{\RR}$ be a
  lattice, and let $C\subseteq N_{\RR}$ be a rational, polyhedral
  strongly convex cone. If $S=C\cap N$, and if $m_1,\ldots,m_p$ are
  the first non-zero integral vectors on the rays of $C$, then there
  is a positive integer $d$ such that for every $m\in S$, $dm$ lies in
  the semigroup $T$ generated by $m_1,\ldots,m_p$.
\end{lemma}

\begin{proof} Consider a smooth fan $\Delta$ that refines the cone
  $C$.  By taking the first non-zero integral vectors on the rays in
  $\Delta$, we get $t$ extra vectors ${m}'_1,\ldots,{m}'_t$. Since each of the cones in
$\Delta$ is smooth, it follows that $S$ is
equal to the semigroup generated by the $m_i$
and the ${m}'_j$.
On the other hand, it is clear that the $m_i$
span the cone
$C\cap N_{\QQ}$ over $\QQ$. Therefore for every
$j\leq t$, we can find
$d_j\in\NN\setminus\{0\}$  such that $d_j{m}'_j$
is in
$T$. Take $d$ to be the least common multiple
of the $d_j$.
\end{proof}

\begin{proof}[Proof of
Proposition~\ref{MDS_algebraic}]
We have already noticed that it is enough
to prove the statement when $V={\rm Spec}(R_0)$ is affine.
Moreover, after taking a refinement of $S$, we may 
assume that the cone $C$ spanned by $S$ is strongly convex.
We use induction on $\dim(S)$. 
If $\dim(S)=1$, then we are done by
Lemma~\ref{dim_one}. 

Suppose now that
$\dim(S)>1$ and that we know the assertion
in smaller dimensions. We use the
construction in the proof of
Lemma~\ref{subcone}. Let $R=R(\ba_{\bullet})$
be the Rees algebra of $\ba_{\bullet}$,
and
let $x_1,\ldots,x_q$ be
homogeneous generators of $R$
as an $R_0$-algebra. We put
$m_i=\deg(x_i)$. Consider the surjective
homomorphism of
$R_0$-algebras
$\Phi : R_0[X_1,\ldots,X_q]\longrightarrow R$,
given by
$\Phi(X_i)=x_i$, and the corresponding semigroup
homomorphism
$\phi : \NN^q\longrightarrow S$ which takes the
$i^{\rm th}$ coordinate vector to $m_i$. Let
$\phi_{\RR}$ be the extension of $\phi$ as a
map $\RR^q\longrightarrow N_{\RR}$.

Consider a smooth fan $\Delta$   refining $C$
such that every
$m_i$  is on a ray of $\Delta$. We apply now
the induction hypothesis for each cone in
$\Delta$ of dimension $\dim(S)-1$ (note that
Lemma~\ref{subcone} ensures the finite
generation of the corresponding
$R_0$-subalgebras). By refining
$\Delta$, we may assume that each face of
dimension $\dim(S)-1$ (as well as its
refinements) satisfies
(\ref{inclusion_in_product}) for a given
positive integer $d$. For example, we take $d$
to be the least common multiple of the positive
integers we get for each face. 
Note that every 
refinement of such
$\Delta$ still satisfies these conditions
(for a possibly different $d$). 
In order to complete
the induction step, it
is enough to show that every  maximal cone
$\sigma\in\Delta$ satisfies
(\ref{inclusion_in_product}) for this $d$.

Let $e_1,\ldots,e_s$ be the generators of
$S_{\sigma}:=\sigma\cap N$ (hence $s=\dim(S)$).
We put
$\widetilde{S}_{\sigma}:=
\phi^{-1}(S_{\sigma})$, and
$\widetilde{\sigma}:=\phi_{\RR}^{-1}(\sigma)$
so that
$\widetilde{S}_{\sigma}=\widetilde{\sigma}\cap\NN^q$.

It is clear that $\widetilde{\sigma}$ is a
rational, polyhedral, strongly convex cone. Now we claim that
every element on ray of $\widetilde{\sigma}$ is mapped by
$\phi_{\RR}$ to the boundary of
$\sigma$. Indeed, suppose that
$w$ is nonzero and lies on a ray of $\widetilde{\sigma}$.
If $w$ is also on a ray of $\RR_+^q$, then $\phi_{\RR}(w)$ is also on a
ray of $\sigma$ by our construction. Otherwise, $w$ is in the interior
of an $r$-dimensional face $F$ of $\RR_{+}^q$, where $2 \le r \le q$.
If $\phi_{\RR}(w)$ is in the interior of $\sigma$, then since 
$\phi_{\RR}$ is continuous and $\sigma$ is of maximal dimension, we can
find an open convex neighborhood $V$ of $w$ in $F$, such that $\phi_{\RR}(V)$
is contained in the interior of $\sigma$. But this contradicts 
the fact that $w$ 
lies on a ray of $\widetilde{\sigma}$. We conclude that $\phi_{\RR}(w)$
is in the boundary of $\sigma$.

We apply now Lemma~\ref{multiple_in_cones} to find $d'$
such that every element in $(d'\cdot \NN)^q\cap
\widetilde{S}_{\sigma}$ is in the semigroup
generated by the first integral points on the
rays of 
$\widetilde{\sigma}$.

Suppose that $f\in \ba_{\sum_ip_ie_i}$, with
$p_i\in d \cdot \NN$. Since
$\Phi$ is surjective, we can write
$f=\sum_{\alpha}c_{\alpha}f_{\alpha}$, where
$c_{\alpha}\in R_0$ and each $f_{\alpha}$ is of
the form  $\Phi(X^u)$, with 
$u\in \phi^{-1}(\deg(f))\subseteq
\widetilde{S}_{\sigma}$. Since
$d'u$ lies in the semigroup generated by the
first integral points on the rays of
$\widetilde{S}_{\sigma}$, it follows that we
can write $f_{\alpha}^{d'}=\prod_ig_i$, where
each $g_i$ is homogeneous, and
$\deg(g_i)=\sum_j\theta_{ij}e_j$  lies in the
boundary of $\sigma$. It follows from the
induction hypothesis that 
$$g_i^d\in\overline{\prod_j
\ba_{de_j}^{\theta_{ij}}},$$ so that
$$f_{\alpha}^{dd'}\in\overline{\prod_j\ba_{de_j}^{d'p_j}}.$$
Since
$d\vert p_j$ for every $j$, we deduce
$$f_{\alpha}\in\overline{\prod_j\ba_{de_j}^{p_j/d}}.$$
Since
$f=\sum_{\alpha}c_{\alpha} f_{\alpha}$, this
implies that $\fra_{\sum_ip_ie_i}\subseteq\overline{\prod_j\fra_{de_j}^{p_j/d}}$.
As we clearly have the inclusion $\prod_j\fra_{de_j}^{p_j/d}\subseteq
\fra_{\sum_ip_ie_i}$, this completes the proof.
\end{proof}

We apply now Proposition~\ref{MDS_algebraic} to prove Theorem
\ref{MDS_local}.

\begin{proof}[Proof of
    Theorem~\ref{MDS_local}] Consider the set  $C$ consisting of
  those $m=(m_i)\in\QQ^r$ such that 
$$h^0(X,{\mathcal
    O}_X(pm_1D_1+\ldots+pm_rD_r))\neq 0$$
for some positive integer
  $p$ with $pm_i\in\ZZ$ for all $i$. It is clear that $C$ is the set
  of points in $\QQ^r$ of a rational convex cone. If we take a finite
  set of homogeneous generators of ${\rm Cox}(D_1,\ldots,D_r)$ as a
  $\CC$-algebra, then their degrees span $C$, so $C$ is
  polyhedral. Denote by $\overline{C}$ the closure of $C$ in
  $\RR^r$.

We have the following
  inclusions 
$$\overline{\phi_{\QQ}^{-1}(\Biggg(X)_{\QQ})}\subseteq\overline{C}
\subseteq\phi_{\RR}^{-1}(\overline{\Biggg(X)_{\RR}}).$$
Since we have assumed that some linear combination of the $D_i$ 
is big, we deduce that the above inclusions are equalities, so
$\overline{B}= \overline{C}$, and therefore it is polyhedral.

We consider now the $\ZZ^r$-graded system
$\frb_{\bullet}=(\frb_m)_{m\in\ZZ^r}$, where $\frb_m$ defines the base
locus of 
$|m_1D_1+\ldots+m_rD_r|$. Our hypothesis implies that this is a
finitely generated system, so we can find a fan $\Delta$ refining
$\ZZ^r$ as in Proposition~\ref{MDS_algebraic}. If $v$
is a discrete valuation of the function field of $X$,
then we define as in \S 2
$$\widetilde{v} : \overline{B}\cap\QQ^r\longrightarrow\RR_+$$
$$\widetilde{v}(m)=
\lim_{p\to\infty}\frac{v(\linser{pm_1D_1+\ldots+pm_rD_r})}{p},$$
where the limit is over those $p$ which are divisible enough. 
Since the valuation of an ideal is equal to that of its integral
closure, it follows from (\ref{inclusion_in_product}) that this
function is linear on each cone in $\Delta$. It 
follows that $\widetilde{v}$ can be uniquely extended by
continuity to $\overline{B}$ (and the extension is again piecewise
linear). Moreover, it is clear from definition that
$\widetilde{v}$ agrees with the pull-back of $v(\parallel\cdot\parallel)$ on $B$.
\end{proof}

\begin{remark}
If $X$ is smooth, similar considerations apply to the functions
${\rm Arn}_Z$ and $\ee_Z$ introduced at the end of \S2. 
Note however that the function 
$m=(m_i)\in\ZZ^r\longrightarrow\Arn_Z(\ba_1^{m_1}\cdot\ldots\cdot\ba_r^{m_r})$
is
not necessarily linear. It is however piecewise linear (it is linear
on a fan refinement 
which does not depend on $Z$, but only on the log resolution of the
ideals $\ba_1,\ldots,\ba_r$). 
Therefore we get our conclusion after passing to a suitable refinement
of $\Delta$.

In the case of $\ee_Z$, it follows from
(\ref{inclusion_in_product}) that the set of those $m\in\QQ^r$ such
that 
$Z$ is not properly contained in an irreducible component of
$\BBasym(m_1D_1+\ldots+m_rD_r)$ 
is the set of rational points in a union of cones in $\Delta$. For
such $m$ we define $\widetilde{\ee_Z}(m)$ in the obvious way, and
(\ref{inclusion_in_product}) implies that $\widetilde{\ee_Z}$ is
polynomial of degree $d$ on each of these cones.  
\end{remark}

The case of varieties with finitely generated linear series, which was
stated in the Introduction, follows easily from Proposition~\ref{MDS_algebraic}.

\begin{proof}[Proof of
    Theorem~\ref{Mori.Dream.Space.Polyhedral.Intro.Thm}]
Take divisors $D_1,\ldots,D_r$ as in Definition~\ref{fg_varieties}. 
If we consider the corresponding map $\phi :
\ZZ^r\longrightarrow
N^1(X)$, then $\phi_{\RR}$ is an
isomorphism. All the assertions now follow from
Theorem~\ref{MDS_local}. 
\end{proof}

\begin{remark} In the context of
  Theorem~\ref{Mori.Dream.Space.Polyhedral.Intro.Thm}, note that if
  $L$ is a line bundle whose class $\alpha$ is on the boundary of
  $\overline{{\rm Eff}}(X)_{\RR}$, then 
it is not clear that $v(\parallel\alpha\parallel)$ (or the other functions) 
can be defined in terms of linear series of multiples of
  $L$. On the other hand, it follows from the proof of
  Theorem~\ref{MDS_local} 
that there does exist \emph{some} line bundle $M$ numerically equivalent to $L$ such
  that $v(\parallel\alpha\parallel)$ can be defined using linear
  series of multiples of $M$.
\end{remark}

We conclude with another application of
Proposition~\ref{MDS_algebraic} to the study of the volume
function. We fix a smooth $n$-dimensional variety $X$. 
Recall that if $L\in\Pic(X)$, then the volume of $L$ is given by 
$$\vol(L):=\limsup_{m\to\infty}\frac{n!\cdot h^0(X,L^m)}{m^n}.$$
This induces a continuous function on $N^1(X)_{\RR}$ such 
that $\vol(mL)=m^n\cdot\vol(L)$ and such that $\vol(L)>0$ 
if and only if $L$ is big. For a detailed study of the volume function
we refer to \cite[Chapter 2]{positivity}. 

We will need the following formula for the volume of a line bundle
which is a consequence of 
Fujita's Approximation Theorem (see \cite{del} or \cite[Chapter
  11]{positivity}). If $L$ is a line bundle with $Bs(L)$ defined by
$\frb\neq\OO_X$, and if $\pi : X'\longrightarrow X$ is a projective,
birational morphism, with $X'$ smooth and such that
$\pi^{-1}(\frb)=\OO_{X'}(-F)$ is an invertible ideal, then we put
$(L^{[n]}):=(M^n)$, 
where $M=\pi^*L-F$. If $\frb=\OO_X$, then we put $(L^{[n]})=0$. With
this notation, we have
\begin{equation}\label{formula_for_volume}
\vol(L)=\sup_{m\in\NN}\frac{((mL)^{[n]})}{m^n}.
\end{equation}
Note that in the above definition of $(L^{[n]})$ we may replace the
ideal $\frb$ by its integral closure.

\begin{proposition}\label{volume}
If $X$ has finitely generated linear series then the closed cone 
$\overline{\BBig(X)}_{\RR}$ has a fan refinement $\Delta$ such that 
the volume function is piecewise polynomial with respect to this fan.
\end{proposition}

\begin{proof} The proof is analogous to the proof of
  Theorem~\ref{Mori.Dream.Space.Polyhedral.Intro.Thm}. In fact, we use
  the same fan refinement. 
By Proposition~\ref{MDS_algebraic}, it is enough to prove the
  following assertion: suppose that $L_1,\ldots,L_r$ are line bundles
  on $X$ 
whose classes are linearly independent in $N^1(X)_{\RR}$, 
and let us denote by $\fra_p$ the base ideal of $\sum_ip_iL_i$ for
  $p\in\NN^r$; if there is $d\geq 1$ such that 
\begin{equation}
\overline{\fra_{dp}}=\overline{\prod_i\fra_{pe_i}^{p_i}}
\tag{$\dagger$}
\end{equation}
for all $p\in\NN^n$, then the volume function is polynomial on the
cone spanned by the classes of $L_1,\ldots,L_r$. 

It is clear that it is enough to show that the map
$$p\longrightarrow\vol\Big(\sum_{i=1}^rdp_iL_i\Big)$$
is a polynomial function of degree $n$ for $p\in\NN^r$.
Let $\pi : X'\longrightarrow X$ be a projective birational morphism, 
with $X'$ smooth and such that $\pi^{-1}(\fra_{de_i})=\OO(-F_i)$ 
are invertible for all $i$. If $M_i=\pi^*(dL_i)-F_i$, then it follows
from $(\dagger)$ that for every $p\in\NN^r$ we have
$$\Big(\Big(\sum_idp_iL_i\Big)^{[n]}\Big)=\Big(\Big(\sum_ip_iM_i\Big)^n\Big).$$
Together with (\ref{formula_for_volume}), this implies 
$$\vol\Big(\sum_idp_iL_i\Big)=\Big(\Big(\sum_ip_iM_i\Big)^n\Big),$$
which completes the proof.
\end{proof}

\subsection*{Acknowledgements}
We are grateful to Arnaud Beauville, Andr\'{e} Hirschowitz and
especially Karen Smith for helpful discussions.  
Research of the first, second, third and fifth authors was partially supported by the NSF under grants DMS 0200278, DMS 0139713, DMS 0500127, and DMS 0200150. The third author was also partially supported by a Clay Mathematics Institute Fellowship.

\end{document}